\documentclass[12pt]{article}

\usepackage[dvipdfm, bookmarksopen=true]{hyperref}
\usepackage[centertags]{amsmath}
\usepackage{amsfonts}
\usepackage{amssymb}
\usepackage{amsthm}

\theoremstyle{plain}
\newtheorem{thm}{Theorem}[section]
\newtheorem{cor}[thm]{Corollary}
\newtheorem{lem}[thm]{Lemma}
\newtheorem{prop}[thm]{Proposition}

\theoremstyle{definition}
\newtheorem{defn}[thm]{Definition}

\newtheorem{rem}[thm]{Remark}
\theoremstyle{remark}
\numberwithin{equation}{section}

\newcommand{\beast}{\begin{eqnarray*}}
\newcommand{\eeast}{\end{eqnarray*}}

\title{Combinatorial representations of Coxeter groups over a field of two elements\thanks{Research partially supported by the NSC grant 96-2628-M-009-015 of
Taiwan R.O.C..}}

\author{Hau-wen Huang\footnote{ Email address:
poker80@msn.com (Hau-wen Huang).} \and Chih-wen
Weng\footnote{$~^\dag$Department of Applied Mathematics National
Chiao Tung University 1001 Ta Hsueh Road Hsinchu, Taiwan 300,
R.O.C..}}

\date{April 14, 2008}

\begin{document}
\maketitle

\begin{abstract}
Let $W$ denote a simply-laced Coxeter group with $n$ generators.
We construct an $n$-dimensional representation $\phi$ of $W$ over
the finite field $F_2$ of two elements. The action of $\phi(W)$ on
$F_2^n$ by left multiplication is corresponding to a combinatorial
structure extracted and generalized from Vogan diagrams.  In each
case W of types A, D and E, we determine the orbits of $F_2^n$
under the action of $\phi(W)$, and find that the kernel of $\phi$
is the center $Z(W)$ of $W.$
\end{abstract}


\section{Introduction}\label{s1}

A {\it simply-laced Coxeter group} is a group $W_S(m)$ with a
finite set of generators $S\subseteq W_S(m)$ subject only to
relations
$$(ss')^{m(s, s')}=1,$$
where $m(s, s)=1$ and $m(s, s')=m(s', s)\in\{2, 3\}$ for
$s\not=s'$ in $S$. When $m$ is specified, we write $W_S$ for
$W_S(m)$, and if both $S$ and $m$ are specified, we write $W$ for
$W_S(m)$. A {\it Coxeter graph} $S$ represents a simply-laced
Coxeter group $W_S(m)$, and vice versa. The vertex set of this
graph is $S$, and there is an edge joining two vertices $s$ and
$s'$ whenever $m(s, s')=3.$ There are Coxter groups which are not
simple-laced. In this article we always assume simple-laced
property in a Coxter group to make the corresponding Coxeter graph
$S$ a simple graph.
\bigskip

We shall investigate a kind of {\it flipping puzzle}, which is
also studied in \cite{xw:07,wu:06}, associated with a given
Coxeter graph $S$. The {\it configuration} of the flipping puzzle
is $S$, together with an {\it assignment} of a unique state, white
or black,  on each vertex of $S$. A {\it move} in the puzzle is to
select a vertex $s$ which has black state, and then flip the state
of each neighbor of $s$. When $S$ is one of the Dynkin diagrams
described in Figure 1, the configuration above is essentially a
{\it Vogan diagram with identity involution}, which was first
defined in \cite{k:96}, in a more general way,  as a combinatorial
object representing the real form of the corresponding complex
simple Lie algebra and a system of choices. See also \cite{pb:00,
pb:02, mkc:04}.
\bigskip

~~

\setlength{\unitlength}{1mm}
\begingroup\makeatletter\ifx\SetFigFont\undefined
\def\x#1#2#3#4#5#6#7\relax{\def\x{#1#2#3#4#5#6}}%
\expandafter\x\fmtname xxxxxx\relax \def\y{splain}
\gdef\SetFigFont#1#2#3{%
  \ifnum #1<17\tiny\else \ifnum #1<20\small\else
  \ifnum #1<24\normalsize\else \ifnum #1<29\large\else
  \ifnum #1<34\Large\else \ifnum #1<41\LARGE\else
     \huge\fi\fi\fi\fi\fi\fi
  \csname #3\endcsname}%
\else \gdef\SetFigFont#1#2#3{\begingroup
  \count@#1\relax \ifnum 25<\count@\count@25\fi
  \def\x{\endgroup\@setsize\SetFigFont{#2pt}}%
  \expandafter\x
    \csname \romannumeral\the\count@ pt\expandafter\endcsname
    \csname @\romannumeral\the\count@ pt\endcsname
  \csname #3\endcsname}%
\fi\endgroup

\begin{picture}(100, 90)
\put(0,89){\footnotesize{$A_n(n\geq1)$}}
\put(30,90){\circle{1.5}}
\put(38,90){\circle{1.5}}
\put(46,90){\circle{1.5}}
\put(50,90){\circle*{.8}}
\put(54,90){\circle*{.8}}
\put(58,90){\circle*{.8}}
\put(62,90){\circle{1.5}}
\put(70,90){\circle{1.5}}
\put(78,90){\circle{1.5}}
\put(30.75,90){\line( 1,0){6.5}}
\put(38.75,90){\line( 1, 0){6.5}}
\put(62.75,90){\line( 1,0){6.5}}
\put(70.75,90){\line( 1, 0){6.5}}
\put(29,87){\scriptsize{$s_n$}}
\put(36,87){\scriptsize{$s_{n-1}$}}
\put(45,87){\scriptsize{$s_{n-2}$}}
\put(61,87){\scriptsize{$s_3$}}
\put(69,87){\scriptsize{$s_2$}}
\put(77,87){\scriptsize{$s_1$}}

\put(0,74){\footnotesize{$D_n(n\geq 4)$}}
\put(30,70){\circle{1.5}}
\put(30,80){\circle{1.5}}
\put(38,75){\circle{1.5}}
\put(46,75){\circle{1.5}}
\put(50,75){\circle*{.8}}
\put(54,75){\circle*{.8}}
\put(58,75){\circle*{.8}}
\put(62,75){\circle{1.5}}
\put(70,75){\circle{1.5}}
\put(78,75){\circle{1.5}}
\put(30.6,70.2){\line( 5, 3){6.9}}
\put(30.6,79.8){\line(5,-3){6.9}}
\put(38.75,75){\line( 1, 0){6.5}}
\put(62.75,75){\line(1, 0){6.5}}
\put(70.75,75){\line( 1, 0){6.5}}
\put(29,67){\scriptsize{$s_{n-1}$}}
\put(28,77){\scriptsize{$s_{n}$}}
\put(37,72){\scriptsize{$s_{n-2}$}}
\put(45,72){\scriptsize{$s_{n-3}$}}
\put(61,72){\scriptsize{$s_3$}}
\put(69,72){\scriptsize{$s_2$}}
\put(77,72){\scriptsize{$s_1$}}

\put(0,55){\footnotesize{$E_6$}}
\put(30,55){\circle{1.5}}
\put(38,55){\circle{1.5}}
\put(46,55){\circle{1.5}}
\put(54,55){\circle{1.5}}
\put(62,55){\circle{1.5}}
\put(46,63){\circle{1.5}}
\put(30.75,55){\line( 1, 0){6.5}}
\put(38.75,55){\line( 1, 0){6.5}}
\put(46.75,55){\line( 1, 0){6.5}}
\put(54.75,55){\line( 1, 0){6.5}}
\put(45.75,55.75){\line( 0, 1){6.5}}
\put(30,52){\scriptsize{$s_5$}}
\put(38,52){\scriptsize{$s_4$}}
\put(46,52){\scriptsize{$s_3$}}
\put(54,52){\scriptsize{$s_2$}}
\put(62,52){\scriptsize{$s_1$}}
\put(46,65){\scriptsize{$s_6$}}

\put(0,35){\footnotesize{$E_7$}}
\put(30,35){\circle{1.5}}
\put(38,35){\circle{1.5}}
\put(46,35){\circle{1.5}}
\put(54,35){\circle{1.5}}
\put(62,35){\circle{1.5}}
\put(46,43){\circle{1.5}}
\put(70,35){\circle{1.5}}
\put(30.75,35){\line( 1,0){6.5}}
\put(38.75,35){\line( 1, 0){6.5}}
\put(46.75,35){\line( 1,0){6.5}}
\put(54.75,35){\line( 1, 0){6.5}}
\put(46,35.75){\line( 0,1){6.5}}
\put(62.75,35){\line( 1, 0){6.5}}
\put(30,32){\scriptsize{$s_6$}}
\put(38,32){\scriptsize{$s_5$}}
\put(46,32){\scriptsize{$s_4$}}
\put(54,32){\scriptsize{$s_3$}}
\put(62,32){\scriptsize{$s_2$}}
\put(46,45){\scriptsize{$s_7$}}
\put(70,32){\scriptsize{$s_1$}}

\put(0,15){\footnotesize{$E_8$}}
\put(30,15){\circle{1.5}}
\put(38,15){\circle{1.5}}
\put(46,15){\circle{1.5}}
\put(54,15){\circle{1.5}}
\put(62,15){\circle{1.5}}
\put(46,23){\circle{1.5}}
\put(70,15){\circle{1.5}}
\put(78,15){\circle{1.5}}
\put(30.75,15){\line( 1, 0){6.5}}
\put(38.75,15){\line( 1, 0){6.5}}
\put(46.75,15){\line( 1, 0){6.5}}
\put(54.75,15){\line( 1, 0){6.5}}
\put(46,15.75){\line( 0, 1){6.5}}
\put(62.75,15){\line( 1, 0){6.5}}
\put(70.75,15){\line( 1, 0){6.5}}
\put(30,12){\scriptsize{$s_7$}}
\put(38,12){\scriptsize{$s_6$}}
\put(46,12){\scriptsize{$s_5$}}
\put(54,12){\scriptsize{$s_4$}}
\put(62,12){\scriptsize{$s_3$}}
\put(46,25){\scriptsize{$s_8$}}
\put(70,12){\scriptsize{$s_2$}}
\put(78,12){\scriptsize{$s_1$}}

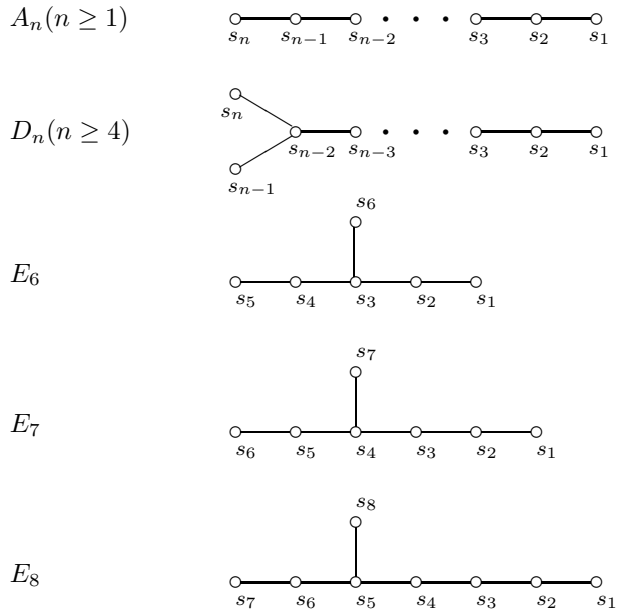
\put(0,3){{\bf Figure 1:} Simply-laced Dynkin diagrams.}
\end{picture}
\bigskip

We fix a simply-laced Coxeter group $W$ and its Coxeter graph $S,$
where $|S|=n.$ Let $F_2$ denote the finite field of two elements
$0$ and $1$. In Section~\ref{s2}, we use the column vector set
$F_2^S=F_2^n$ to describe the set of configurations in the
flipping puzzle associated with $S$ by setting that $\ell_s=1$ iff
the configuration $\ell\in F_2^n$ has black state in the vertex
$s.$ For each vertex $s\in S$, we find a way to associate the move
with selecting vertex $s$ as an $n\times n$ invertible matrix
$\mathbf{s}$ over $F_2.$ This $\mathbf{s}$ acts on a configuration
$\ell\in F_2^n$ by left multiplication to become a new
configuration $\mathbf{s}\ell$ which has the desired property as
stated in the definition of the flipping puzzle when $\ell_s=1.$
Unlike in the definition, our move $\mathbf{s}$ does not select
configuration $\ell$, but if a configuration has white state in
$s$, it makes no effect; i.e. if $\ell_s=0$ then
$\mathbf{s}\ell=\ell.$
\bigskip

Let ${\rm GL}_n(F_2)$ denote the set of $n\times n$ invertible
matrices over $F_2$ and let $\mathbf{W}$ denote the subgroup of
${\rm GL}_n(F_2)$ generated by the moves $\mathbf{s}$ for $s\in
S.$ We refer $\mathbf{W}$ to a {\it flipping group} of $S.$ In
Section~\ref{s3}, we find that the canonical map
$\phi:W\rightarrow GL_n(F_2),$ lifted from $\phi(s)=\mathbf{s}$
for $s\in S$, is a homomorphism with $\phi(W)=\mathbf{W}.$ Due to
its origination, we refer such a map to the {\it Vogan
representation} of $W.$ Then we find that the flipping group
$\mathbf{W}$ has trivial center in Section~\ref{s4}. In Sections
\ref{s5}, \ref{s6} and \ref{s7}, we assume $W$ to be $A_n$, $D_n$
and $E_n$ respectively. By using the finiteness of $W$, we can
determine the size of the corresponding flipping group
$\mathbf{W}$. We find that the kernel of the Vogan representation
of $W$ is the center $Z(W)$ of $W$ when $W$ is finite.
\bigskip

In the flipping puzzle on a Coxeter graph $S$, two configurations
are said to be {\it equivalent} if one can be obtained from the
other by a sequence of moves. Let $\mathcal{P}$ denote the
partition of configurations (i.e. $F_2^n$) according to the above
equivalent relation. As a byproduct of our work, we solve the
flipping puzzle associated with $S$ when $S$ is each of $A_n$,
$D_n$ and $E_n$ by determining $\mathcal{P}$. Note that when $S$
is a {\it tree}, a generalization of Dynkin diagrams, some partial
results on $\mathcal{P}$ are obtained in \cite{wu:06} and
\cite{xw:07}.
\bigskip

\section{Flipping groups}\label{s2}

Throughout this article,  $W$ will be a simply-laced Coxeter group
with corresponding Coxter graph $S$ of $n$ elements and edge set
$R=\{ss'~|~m(s, s')=3\}$. We shall construct a matrix group
associated with the flipping puzzle on the Coxeter graph $S$. Let
${\rm Mat}_n(F_2)$ denote the set of $n\times n$ matrices over
$F_2$ with rows and columns indexed by $S.$ Let $F_2^n$ denote the
set of $n$-dimensional column vectors over $F_2$ indexed by $S$.
For $s\in S$, let $\widetilde{s}$ denote the characteristic vector
of $s$ in $F_2^n$; that is $\widetilde{s}=(0, 0, \ldots, 0, 1, 0,
\ldots, 0)^t,$ where $1$ is in the position corresponding to $s.$
\bigskip

\begin{defn}\label{d2.1}
For $s\in S,$ we associate a matrix $\mathbf{s}\in {\rm
Mat}_n(F_2)$, denoted by the bold type of $s$, as
$$\mathbf{s}_{uv}=\left\{
\begin{array}{ll}
    1, & \hbox{if $u=v,$ or $v=s$ and $uv\in R$;} \\
0, & \hbox{else,} \\\end{array} \right.$$ where $u,v\in S$.
\end{defn}

The following is a reformulating of Definition~\ref{d2.1}.

\begin{lem}\label{l2.2} For $s, v\in S,$
$$\mathbf{s}\widetilde{v}=\left\{%
\begin{array}{ll}
  \widetilde{v}  , & \hbox{if $v\not= s$;} \\ %
\widetilde{v}+\sum\limits_{uv\in R}\widetilde{u} & \hbox{if $v=s$.} \\
\end{array}%
\right.
$$\hfill $\Box$
\end{lem}


The flipping puzzle associated with $S$, which is described in the
introduction, is now restated as follows. A configuration is
simply a vector $\ell\in F_2^n,$ where $\ell_s=1$ (resp.
$\ell_s=0$) means that the vertex $s\in S$ has black state (resp.
white state). In this setting, if $\ell_s=1$ then $\mathbf{s}\ell$
is the new configuration after the move to select the vertex $s$.
Note that if $\ell_s=0$, we have $\mathbf{s}\ell=\ell$ from
Lemma~\ref{l2.2}, so we can view the action of $\mathbf{s}$ on
$\ell$ as a {\it feigning move} on $\ell$ which is not originally
defined as a move in the flipping puzzle.  The following lemma is
immediate from this combinatorial realization.

\begin{lem}\label{l2.3} For $s \in S,$ $\mathbf{s}$ is an involution;
that is  ${\mathbf{s}}^2=I$, the identity matrix. \hfill $\Box$
\end{lem}

From Lemma~\ref{l2.3}, $\mathbf{s}$ is invertible, so we can give
the following definition.

\begin{defn}\label{d2.4}
Let $\mathbf{W}$ denote the subgroup of  ${\rm GL}_n(F_2)$ generated
by the set $\{\mathbf{s}~|~s\in S\}$.  $\mathbf{W}$ is referring to
the {\it flipping group} of $S$.
\end{defn}

\section{Coxeter groups and their combinatorial
representations}\label{s3}

Let $W$ denote a simply-laced Coxeter group. Recall that an {\it
$n$-dimensional representation} of $W$ over  $F_2$ is a
homomorphism of $W$ into ${\rm GL}_n(F_2).$ It is notorious
difficult in the study of groups only defined by generators and
relations. Hence the representation theory of Coxeter groups plays
an important role in the study. In \cite[Section~5.3]{h:90},
Humphreys gives "geometric representations" of Coxeter groups and
use these representations to show that the finite Coxeter groups
are essentially those associated with Dynkin diagrams. In this
section we shall show that the flipping groups defined in the last
section give "combinatorial representations" of simply-laced
Coxeter groups. First we need a lemma.

\begin{lem}\label{l3.2}
Let $W$ denote a simply-laced Coxeter group with Coxeter graph
$S.$ For $s\in S$, set $E_s\in {\rm Mat}_n(F_2)$ by
\begin{equation}\label{e3.1}E_s\widetilde{v}=\left\{%
\begin{array}{ll}
    0, & \hbox{if $v\not=s$;} \\
\sum\limits_{uv\in R}\widetilde{u}, & \hbox{if $v=s$} \\
\end{array}%
\right. \qquad {\rm for}~v\in S.\end{equation} Then with referring
to the notation in Definition~\ref{d2.1}, the following (i)-(iii)
hold.
\begin{enumerate}
\item[(i)] $\mathbf{s}=I+E_s$ for $s\in S$ \item[(ii)]
$E_{s'}E_s=0$, if $s's\notin R$. \item[(iii)] If $s_is_{i-1}\in R$
for $i=1,2,\ldots,t$, then
$$
 E_{s_t}E_{s_{t-1}}\cdots E_{s_0}=\left\{%
\begin{array}{ll}
E_{s_0}, & \hbox{if $s_t=s_0$;} \\
E_{s_t}E_{s_0}, & \hbox{if $s_ts_0\in R$.}\\
\end{array}%
\right.
$$
\end{enumerate}
\end{lem}
\begin{proof}
(i) is immediate from Lemma~\ref{l2.2}. Note that
$E_{s'}E_s\widetilde{v}=0$  by (\ref{e3.1}) for any $v, s, s'\in S$
with $s's\not\in R,$ and hence we have (ii). (iii) follows from the
same reason as in (ii) by applying the product of matrices in either
side of the equation to $\widetilde{v}$ and obtaining the desired
equality in each case.
\end{proof}

\begin{thm}\label{t3.2} Let $W$ denote a simply-laced Coxeter group with Coxeter graph $S$.
Let $\mathbf{W}$ denote the flipping group of $S$. Then there
exists a surjective homomorphism $\phi:W\rightarrow \mathbf{W}$
such that $\phi(s)=\mathbf{s}$ for $s\in S.$ In particular, $\phi$
is a representation of $W$ over $F_2.$
\end{thm}
\begin{proof}
We have seen $\mathbf{s}^2=I$ for $s\in S.$ It remains to show
$(\mathbf{s}\mathbf{s}')^2=I$ if $s\neq s'$ and $ss'\not\in R,$ and
to show $(\mathbf{s}\mathbf{s}')^3=I$ if $ss'\in R.$ For $s, s'\in
S$,
\begin{eqnarray*}
\mathbf{s}\mathbf{s}'&=&(I+E_s)(I+E_{s'}) \\
    &=&I+E_s+E_{s'}+E_sE_{s'}
\end{eqnarray*}
by Lemma~\ref{l3.2}(i). In the case $s\neq s'$ and $ss'\not\in R$,
\begin{eqnarray*}
(\mathbf{s}\mathbf{s}')^2&=&(I+E_s+E_{s'})(I+E_s+E_{s'}) \\
    &=&I+2E_s+2E_{s'}\\
    &=&I
\end{eqnarray*}
by Lemma~\ref{l3.2}(ii). In the case $ss'\in R$,
\begin{eqnarray*}
(\mathbf{s}\mathbf{s}')^2&=&(I+E_s+E_{s'}+E_sE_{s'})(I+E_s+E_{s'}+E_sE_{s'}) \\
    &=&I+3E_s+3E_{s'}+4E_sE_{s'}+E_{s'}E_s \\
    &=&I+E_s+E_{s'}+E_{s'}E_s
\end{eqnarray*}
and
\begin{eqnarray*}
(\mathbf{s}\mathbf{s}')^3&=&(\mathbf{s}\mathbf{s}')^2(\mathbf{s}\mathbf{s}')\\
    &=&(I+E_s+E_{s'}+E_{s'}E_s)(I+E_s+E_{s'}+E_sE_{s'})\\
    &=&I+2E_s+4E_{s'}+2E_sE_{s'}+2E_{s'}E_s\\
    &=&I
\end{eqnarray*}
by Lemma~\ref{l3.2}(iii).
\end{proof}

\begin{defn}\label{d3.4}
The representation $\phi$ defined in Theorem~\ref{t3.2} is called
the {\it Vogan representation} of $W$.
\end{defn}


Suppose $J\subseteq S$. Let $\mathbf{W}_J$ denote the subgroup of
$\mathbf{W}$ generated by the set $\{\mathbf{s}~|~s\in J\}$ and
$W_J$ denote simply-laced Coxeter group with the set $J$ of
generators with the function $m\upharpoonright J\times J$,
 the restriction of $m$ to $J\times J.$
Note that $W_J$ is isomorphic to the subgroup of $W$ generated by
the set $\{s~|~s\in J\}$ \cite[Section~5.5]{h:90}. Hence we use
the same symbol $W_J$ to express these two isomorphic groups. It
makes no confused if the place that $W_J$ appears is also
considered. For example, the first $W_J$ in (iii) of the following
lemma is in the first meaning and the remaining two $W_J$ are in
the second meaning. Note that $\mathbf{W_J},$ which is different
to $\mathbf{W}_J$, is the flipping group on $J.$ Let $G[J]$ denote
the submatrix of $G\in {\rm Mat}_n(F_2)$ with rows and columns
indexed by $J$, and $\mathbf{W}_J[J]:=\{G[J]~|~G\in
\mathbf{W}_J\}$.

\begin{lem}\label{l3.5}
Suppose $J\subseteq S.$ The following (i)-(iii) hold.
\begin{enumerate}
\item[(i)] $\mathbf{W}_J[J]=\mathbf{W_J}$. \item[(ii)] The map
$\psi:\mathbf{W}_J\rightarrow \mathbf{W_J}$, defined by
$\psi(G)=G[J]$ for $G\in \mathbf{W}_J$, is a surjective
homomorphism.
\item[(iii)]
 Let $\phi$ and
$\phi'$ denote the Vogan representations of $W_S$ and $W_J$
respectively. Then $\phi'=\psi\circ \phi\upharpoonright W_J$. In
particular, ${\rm Ker}~\phi\upharpoonright W_J\subseteq{\rm
Ker~}\phi'.$
\end{enumerate}
\end{lem}
\begin{proof}
By Definition~\ref{d2.1}, $\mathbf{s}_{uv}=0$ for $s,u\in J$ and
$v\in
S-J$. By this, 
 each matrix $G\in \mathbf{W}_J$ has the form
$$
G=\left(
  \begin{array}{clr}
   A & 0\\
   B & C
  \end{array}
  \right)
$$
if indices in $J$ are placed in the beginning of rows and columns,
where $A$ is a $|J|\times |J|$ matrix, $B$ is an $(n-|J|)\times
|J|$ matrix, $C$ is an $(n-|J|)\times (n-|J|)$ matrix, and $0$ is
a $|J|\times (n-|J|)$ zero matrix. Then (i) and (ii) follow from
the following matrix product rule in block form:
$$
\left(
  \begin{array}{clr}
   A & 0\\
   B & C
  \end{array}
\right) \left(
  \begin{array}{clr}
   A' & 0\\
   B' & C'
  \end{array}
\right) = \left(
  \begin{array}{clr}
   AA' & 0\\
   BA'+CB' & CC'
  \end{array}
\right).
$$
Since $\psi\circ\phi(s)=\mathbf{s}[J]=\phi'(s)$ by (i) for all
$s\in J$, we see $\phi'=\psi\circ \phi\upharpoonright W_J$, and
this implies ${\rm Ker}~\phi\upharpoonright W_J\subseteq{\rm
Ker~}\phi'.$
\end{proof}

\section{The center of a flipping group}\label{s4}

As we shall see in Proposition~\ref{p5.12} that the Coxeter group
of type $D_n$ has nontrivial center when $n$ is even. In this
section, we show that the center $Z(\mathbf{W})$ of any flipping
group $\mathbf{W}$ of a Coxeter graph $S$ is trivial. Therefore,
the center $Z(W)$ of any Coxeter group $W$ is contained in the
kernel of the Vogan representation of $W.$ Recall that a Coxeter
graph $S$ is {\it disconnected} if there is a partition of
$S=S'\cup S''$ with $S', S''\not=\emptyset$ and there is no edge
$uv\in R$ with $u\in S'$ and $v\in S''.$ In this case the Coxeter
group $W$ is isomorphic to the direct product $W'\times W''$ of
the Coxeter groups $W'=W_{S'}$ and $W''=W_{S''}.$ $S$ is {\it
connected} if $S$ is not disconnected.

\begin{prop}\label{p4.1} Let $W$ denote a simple-laced Coxeter group with
Coxeter graph $S.$ Then the center $Z(\mathbf{W})$ of  the
flipping group $\mathbf{W}$ of $S$ is trivial.
\end{prop}
\begin{proof} It suffices to assume that $S$ is connected with at
least two vertices. Let  $Z$ be an element in the center of
$\mathbf{W}$ and let $u,v$ be two distinct elements in $S$. We
show that $Z_{vu}=0$ to conclude $Z=I.$ Suppose $Z_{vu}=1.$ On the
one hand $\mathbf{v}Z\widetilde{u}\not=Z\widetilde{u}$ since
$Z\widetilde{u}$ has $1$ in the $v$th position. On the other hand,
$\mathbf{v}Z\widetilde{u} =
Z\mathbf{v}\widetilde{u}=Z\widetilde{u}$ since
$\mathbf{v}\widetilde{u}=\widetilde{u}.$ Hence we have a
contradiction.
\end{proof}

From the above Proposition~\ref{p4.1} we immediately have the
following corollary.

\begin{cor}\label{c4.2}
Let $W$ denote a simply-laced Coxeter group. Then the center
$Z(W)$ is contained in the kernel of the Vogan representation of
$W.$ \hfill $\Box$
\end{cor}

\section{Coxeter groups of type $A_n$}\label{s5}

Recall that the Vogan representation  $\phi$ of $W$  is {\it
faithful} whenever $\phi$ is injective. Also $\phi$ is {\it
irreducible} if there is no subspace $V\subseteq F_2^n$, $V\not=0,
F_2^n,$ such that $\phi(W)V\subseteq V.$ For $a\in F_2^n$, the
subset of $F_2^n$ consisting of all elements $Ga$ with $G\in
\phi(W)$ is called the {\it orbit} of $F_2^n$ containing $a$ under
the action of $\phi(W).$
\bigskip

In this section we assume that $W$ is of type $A_n$ with the
Coxeter graph $S$ as shown in Fig. 1, and determine the orbits of
$F_2^n$ under the action of $\phi(W)$. We also show that the
kernel of the Vogan representation $\phi$ of $W$ is the center
$Z(W)$ of $W$ and determine the reducibility of $\phi.$ The
trivial case is given in the following.

\begin{prop}
Let $W$ be a Coxeter group of type $A_1$ with the Vogan
representation $\phi$. Then the orbits of $F_2$ are $\{0\}$,
$\{1\}$ under the action of $\phi(W),$
 ${\rm Ker}~\phi=\{1,s_1\}=W=Z(W)$, and $\phi$ is irreducible.
\end{prop}
\begin{proof}
This follows from that $W=\{1,s_1\}$ and $\phi(W)$ is a trivial
group.
\end{proof}

In the remaining of this section, we always assume $n\geq 2$. Set
\begin{equation}\label{e4.1}
\overline{1}=\widetilde{s}_1,~\overline{i+1}=\mathbf{s_i}\mathbf{s_{i-1}}\cdots
\mathbf{s_1}\overline{1}\quad{\rm for}~ 1\leq i\leq n.
\end{equation}
Note that
\begin{equation}\label{e4.0}
\overline{i}=\widetilde{s}_{i-1}+ \widetilde{s}_i \quad {\rm for}~
2\leq i\leq n,
\end{equation}
and
\begin{equation}\label{e4.2}
\overline{n+1}=\widetilde{s}_n=\overline{1}+\overline{2}+\cdots+\overline{n}.
\end{equation}
Set
$\Delta=\Delta(A_n):=\{\overline{1},\overline{2},\ldots,\overline{n}\}.$
Note that $\Delta$ is a basis of $F_2^n.$ We refer $\Delta$ to a
{\it simple basis} of $F_2^n.$ For $a\in F_2^n$, let $\Delta(a)$
denote the subset of $\Delta$ consisting of all the elements
appeared in the expression of $a$ as a linear combination of
elements in $\Delta.$ The {\it weight} of an element $a\in F_2^n$
is $wt(a):=|\Delta(a)|.$ For example,
$\Delta(\overline{n+1})=\Delta$ and $wt(\overline{n+1})=n.$

\begin{lem}\label{l4.2}
$\mathbf{s_i}\overline{i}=\overline{i+1},$
$\mathbf{s_i}\overline{i+1}=\overline{i}$ and $\mathbf{s_i}$ fixes
other vectors in $ \{\overline{1},$ $\overline{2},$ $\ldots,$
$\overline{n+1}\}-\{\overline{i},\overline{i+1}\}$ for $1\leq i\leq
n.$
\end{lem}
\begin{proof}
This is immediate by applying Lemma~\ref{l2.2}, (\ref{e4.1}) and
(\ref{e4.0}).
\end{proof}

Let $S_{n+1}$ denote the group of permutations on
$\{\overline{1},\overline{2},\ldots ,\overline{n+1}\}.$ By Lemma
\ref{l4.2}, we can give the following definition.

\begin{defn}\label{d4.3}
Let $\alpha~:~\mathbf{W}\rightarrow S_{n+1}$ be the homomorphism
defined by
$$\alpha(G)\overline{j}=G\overline{j}$$ for each $1\leq j\leq n+1$ and $G\in \mathbf{W}$.
\end{defn}

Note that $\alpha(\mathbf{s_i})$ is the transposition
$(\overline{i},\overline{i+1})$ in $S_{n+1}$ for each $1\leq i\leq
n$.

\begin{lem}\label{p4.4}
$\alpha$ is an isomorphism from $\mathbf{W}$ onto $S_{n+1}.$
\end{lem}
\begin{proof}
$\alpha$ is surjective since the transpositions
$\alpha(\mathbf{s_1})$, $\alpha(\mathbf{s_2})$,$\ldots$,
$\alpha(\mathbf{s_{n}})$ generate $S_{n+1}$. Since $\Delta\cup
\{\overline{n+1}\}$ spans $F_2^n$, $\alpha$ is injective.
\end{proof}

The next proposition determines the orbits of $F_2^n$ under the
action of $\mathbf{W}$.

\begin{prop}\label{l4.4} For $0 \leq i\leq \lfloor\frac{n+1}{2}\rfloor,$
$$O_i=\{a\in F_2^n~|~wt(a)=i~{\rm or}~n+1-i\}$$ is an orbit of
$F_2^n$ under the action of $\mathbf{W},$ where $\lfloor t\rfloor$
is the largest integer less than or equal to $t.$
\end{prop}
\begin{proof} Suppose $a\in F_2^n$ with $wt(a)=i.$
Observe that from Lemma~\ref{p4.4} and (\ref{e4.2}),
$$\Delta(Ga)=\left\{
\begin{array}{ll}
    \alpha(G)\Delta(a), & \hbox{if $\overline{n+1}\not\in \alpha(G)\Delta(a)$ ;} \\
\Delta-\alpha(G)\Delta(a), & \hbox{if $\overline{n+1}\in \alpha(G)\Delta(a)$}\\
\end{array}%
\right.
$$
for $G\in \mathbf{W}.$ The proposition follows from this
observation because the subgroup of $\alpha(\mathbf{W})=S_{n+1}$
generated by the transpositions $\alpha(\mathbf{s_1})$,
$\alpha(\mathbf{s_2})$,$\ldots$, $\alpha(\mathbf{s_{n-1}})$ acts
transitively on the fixed size subsets of $\Delta$, and
$\mathbf{s_n}\overline{n}=\overline{1}+\overline{2}+\cdots+\overline{n}$
by Lemma~\ref{l4.2} and (\ref{e4.2}).
\end{proof}


In the following propositions, we study the reducibility of $\phi$
and ${\rm Ker}~\phi.$

\begin{prop}\label{p4.6}
The Vogan representation $\phi$ of $W$ is irreducible if and only if
$n$ is even.
\end{prop}
\begin{proof}
Let $V$ denote a nontrivial proper subspace of $F_2^n$ such that
$\phi(W)V\subseteq V$. Referring to Proposition~\ref{l4.4}, note
that \begin{equation}\label{e5.4} V=\bigcup\limits_{i\in J}O_i
\end{equation} for some proper subset $J\subseteq \{0, 1, \ldots,
\lfloor\frac{n+1}{2}\rfloor\}$ with $J\not=\{0\}.$ Note that the
 set in the right side of (\ref{e5.4}) to be closed under addition is
 when it is the set of even weight vectors, and this occurs if and
 only if $n$ is odd.
\end{proof}

\begin{prop}\label{p4.5}
The Vogan representation $\phi$ of $W$ is faithful. In particular,
${\rm Ker}~\phi=Z(W)$ is the trivial group.
\end{prop}

\begin{proof} The first statement follows from that Proposition~\ref{p4.4} and
$W$ is isomorphic to $S_{n+1}$ \cite[p41]{h:90}. The second
follows from the first and Corollary~\ref{c4.2}.
\end{proof}

\section{Coxeter groups of type $D_n$}\label{s6}

Fix an integer $n\geq 4.$ Let $W$ denote the Coxeter group of type
$D_n$ with the Coxeter graph in Fig. 1. Let $\phi$ denote the
Vogan representation of $W$, and $\mathbf{W}=\phi(W)$ be the
flipping group of $S.$ Set
\begin{align}\label{e5.0}
\begin{split}
&\overline{1}=\widetilde{s}_1,~\overline{i+1}=\mathbf{s_i}\mathbf{s_{i-1}}\cdots
\mathbf{s_1}\overline{1}\quad {\rm for}~ 1\leq i\leq n-1,~{\rm
and~} \overline{n+1}=\widetilde{s}_n.
\end{split}
\end{align}
Note that
\begin{align}\label{e5.1}
\begin{split}
&\overline{i}=\widetilde{s}_{i-1}+\widetilde{s}_{i}~~~~{\rm~for}~~2\leq
i\leq n-2,\\
&\overline{n-1}=\widetilde{s}_{n-2}+\widetilde{s}_{n-1}+\widetilde{s}_n,
\end{split}
\end{align}
and
\begin{equation}\label{e5.2}
\overline{n}=\widetilde{s}_{n-1}+\widetilde{s}_n=\overline{1}+\overline{2}+\cdots+\overline{n-1}.
\end{equation}
Set $\Delta=\Delta(D_{n}):=\{\overline{1},
\overline{2},\ldots,\overline{n-1}, \overline{n+1}\}$ to be the
simple basis of $F_2^n$ in the case of type $D_n$. Set $\Delta(a)$
and $wt(a)$ as before for $a\in F_2^n.$ For example,
$\Delta(\overline{n})=\Delta-\{\overline{n+1}\}$ by (\ref{e5.2}),
and $wt(\overline{n})=n-1.$

\begin{lem}\label{l5.1} The following (i),(ii) hold.
\begin{enumerate}
\item[(i)] For each $1\leq i\leq n-1$,
$\mathbf{s_i}\overline{i}=\overline{i+1},$
$\mathbf{s_i}\overline{i+1}=\overline{i},$ and
$$\mathbf{s_i}\overline{j}=\overline{j}\qquad {\rm for~~} j\in
\{\overline{1},\overline{2},\ldots,\overline{n+1}\}-\{\overline{i},\overline{i+1}\}.$$
\item[(ii)] $\mathbf{s_n}\overline{n-1}=\overline{n},$
$\mathbf{s_n}\overline{n}=\overline{n-1},$
$\mathbf{s_n}\overline{n+1}=\overline{n-1}+\overline{n}+\overline{n+1},$
 and
$$\mathbf{s_n}\overline{j}=\overline{j}\qquad {\rm for~~} j\in\{\overline{1}, \overline{2}, \ldots, \overline{n-2}\}.$$
\end{enumerate}
In particular, $\overline{n+1}\in \Delta(G\overline{n+1})$ and
$G(\{\overline{1}, \overline{2}, \ldots, \overline{n}\})\subseteq
\{\overline{1}, \overline{2}, \ldots, \overline{n}\}$ for all $G\in
\mathbf{W}.$

\end{lem}
\begin{proof}
This  follows immediately from Lemma~\ref{l2.2}, (\ref{e5.0}) and
(\ref{e5.1}).
\end{proof}

Let $S_n$ denote the group of permutations on
$\{\overline{1},\overline{2},\ldots,\overline{n}\}$. By
Lemma~\ref{l5.1}, we can give the following definition.

\begin{defn}\label{d5.2}
Let $\beta~:~\mathbf{W}\rightarrow S_n$ denote the homomorphism
defined by
$$\beta(G)(\overline{j})=G\overline{j}$$
for $1\leq j\leq n$ and $G\in\mathbf{W}.$
\end{defn}

In fact, $\beta$ is surjective since the $n-1$ transpositions
$\beta(\mathbf{s_1}),$ $\beta(\mathbf{s_2}),$ $\ldots,$
$\beta(\mathbf{s_{n-1}})$ generate $S_n.$ Let $Z$ denote the
additive group of $(n-1)$-dimensional subspace of $F_2^n$ spanned by
the set $\{\overline{1},\overline{ 2}, \ldots, \overline{n-1}\}.$
Note that $a\in Z$ iff $\overline{n+1}\not\in \Delta(a)$ for $a\in
F_2^n.$ By Lemma \ref{l5.1} and (\ref{e5.2}), $Z$ is closed under
the left multiplication of matrices in $\mathbf{W}.$

\begin{prop}
The Vogan representation $\phi$ of $W$ is not irreducible. In
particular $\phi(W) Z\subseteq Z.$ \hfill $\Box$
\end{prop}

Hence $Z$ is a disjoint union of orbits of $F_2^n$ under the action
of $\mathbf{W}.$ Note that $F_2^n-Z$ is also a disjoint union of
orbits of $F_2^n$ under the action of $\mathbf{W}$. The following
proposition determines all the orbits of $F_2^n$ under the action of
$\mathbf{W}.$

\begin{prop}\label{p5.3} The following are orbits of $F_2^n$ under the action of $\mathbf{W}.$
\begin{eqnarray*}
O_i&=&\{a\in Z~|~wt(a)=i~{\rm or}~n-i\},\\
\Omega_o&=&\{a\in F_2^n-Z~|~wt(a)\equiv 1~{\rm
or}~n-1~\pmod{2}\},\\
\Omega_e&=&\{a\in F_2^n-Z~|~wt(a)\equiv 0~{\rm or}~n~\pmod{2}\},
\end{eqnarray*}
where $0\leq i\leq \lfloor\frac{n}{2}\rfloor.$ In particular
$\Omega_o=\Omega_e=F_2^n-Z$ is an orbit when $n$ is odd.
\end{prop}
\begin{proof} The proof is similar to the proof of
Proposition~\ref{l4.4}. The reason that $O_i$ is an orbit follows
from two facts: (i) $\beta(\mathbf{s_1})$, $\beta(\mathbf{s_2}),$
$\ldots,$ $\beta(\mathbf{s_{n-2}})$ generate the subgroup $S_{n-1}$
of $S_n$ consisting of  permutations on $\Delta-\{\overline{n+1}\}$
and $S_{n-1}$ acts transitively on fixed size subsets of
$\Delta-\{\overline{n+1}\}$, and (ii)
$$\mathbf{s_{n-1}}\overline{n-1}=\mathbf{s_{n}}\overline{n-1}=\overline{n}=
\overline{1}+\overline{2}+\cdots+\overline{n-1}$$ by Lemma
\ref{l5.1}(i),(ii) and (\ref{e5.2}). The reason that $\Omega_o$ and
$\Omega_e$ are orbits follows from an additional fact that
$$wt(\mathbf{s_n}\overline{n+1})=wt(\overline{n-1}+\overline{n}+\overline{n+1})
=wt(\overline{1}+\overline{2}+\cdots+\overline{n-2}+\overline{n+1})=n-1.$$
\end{proof}


We study the structure of $\mathbf{W}.$

\begin{defn}\label{d5.6}
Let $\gamma~:~\mathbf{W}\rightarrow {\rm Aut}(Z)$ denote the
homomorphism from $\mathbf{W}$ into the group ${\rm Aut}(Z)$ of
automorphisms of $Z$ such that
$$\gamma(G)(u)=Gu$$
for $G\in \mathbf{W}$ and $u\in Z.$
\end{defn}

\begin{lem}\label{l5.7}
There exists a unique homomorphism $\theta: S_n\rightarrow {\rm Aut}(Z)$
such that $\gamma=\theta \circ \beta.$
\end{lem}
\begin{proof} Since $\beta$ is surjective,
it suffices to show that the kernel of $\beta$ is contained in the
kernel of $\gamma.$ Suppose $G\in {\rm Ker}~\beta$. Then
$G\overline{i}=\overline{i}$ for $1\leq i\leq n;$ in particular
$G$ fixes each element in the basis $\Delta-\{\overline{n+1}\}$ of
$Z$. Thus $G\in{\rm Ker}~\gamma.$
\end{proof}

Let $Z\rtimes_{\theta} S_n$ denote the group of {\it external
semidirect product} of $Z$ and $S_n$ with respect to
$\theta$\cite[p.155]{dj:02}; i.e. $Z\rtimes_{\theta} S_n$ is the
set $Z\times S_n$ with  the following product rule:
$$(u,\sigma)(v,\tau)=(u+\theta(\sigma)(v),\sigma\tau),$$
where $u, v\in Z$ and $\sigma, \tau\in S_n.$  Note that
$\overline{n+1}+G\overline{n+1}\in Z$ for any $G\in \mathbf{W}$ by
Lemma \ref{l5.1}.

\begin{defn}\label{d5.8}
Let $\delta~:~\mathbf{W}\rightarrow Z\rtimes_{\theta} S_n$ denote
the map defined by
$$\delta(G)=(\overline{n+1}+G\overline{n+1},\beta(G))$$
for any $G\in \mathbf{W}.$
\end{defn}

\begin{lem}\label{l5.9}
$\delta$ is an injective homomorphism of $\mathbf{W}$ into
$Z\rtimes_{\theta} S_n$.
\end{lem}
\begin{proof} For $G, H\in \mathbf{W},$
\begin{eqnarray*}
\delta(G)\delta(H)&=&(\overline{n+1}+G\overline{n+1},
\beta(G))(\overline{n+1}+H\overline{n+1}, \beta(H))\\
&=& (\overline{n+1}+G\overline{n+1}+\theta(\beta(G))(\overline{n+1}+H\overline{n+1}),~\beta(G)\beta(H))\\
&=& (\overline{n+1}+G\overline{n+1}+G(\overline{n+1}+H\overline{n+1}),~\beta(G)\beta(H))\\
&=& (\overline{n+1}+GH\overline{n+1},~\beta(GH))\\
&=&\delta(GH).
\end{eqnarray*}
This shows that $\delta$ is a homomorphism. $\delta$ is injective
since if $\overline{n+1}+G\overline{n+1}=0$ and $G\in {\rm
Ker}~\beta$ , then $G$ fixes all vectors in $\Delta$, so $G$ is the
identity matrix.
\end{proof}

Note that $Z=\overline{n+1}+\Omega_o$ if $n$ is odd, and
$Z=(\overline{n+1}+\Omega_o)\cup (\overline{n+1}+\Omega_e)$ if $n$
is even.

\begin{lem}\label{p5.10} $\delta(\mathbf{W})=(\overline{n+1}+\Omega_o)\rtimes_{\theta} S_n.$
 In particular $\delta(\mathbf{W})=Z\rtimes_{\theta} S_n$ if $n$ is
odd; $\delta(\mathbf{W})$ has index $2$ in $Z\rtimes_{\theta} S_n$
if $n$ is even.
\end{lem}
\begin{proof} Note that $\delta(\mathbf{s_1}), \delta(\mathbf{s_2}), \ldots,
\delta(\mathbf{s_{n-2}}), \delta(\mathbf{s_{n-1}})$ generate
$0\rtimes_{\theta} S_n.$ Since $\Omega_o$ is an orbit containing
$\overline{n+1},$ we have
$\delta(\mathbf{W})=(\overline{n+1}+\Omega_o)\rtimes_{\theta} S_n.$
The second part follows from Proposition~\ref{p5.3}.
\end{proof}

\begin{prop}\label{p5.12}
The Vogan representation $\phi$ of $W$ is faithful when $n$ is odd;
${\rm Ker}~\phi$ has order 2 when $n$ is even. Moreover, ${\rm
Ker}~\phi$ is the center $Z(W)$ of $W.$
\end{prop}
\begin{proof}
Note that $W$ is isomorphic to the semidirect product $Z\rtimes
S_n$ of $Z$ and $S_n$ \cite[p.42]{h:90}. By Lemma~\ref{p5.10},
$\phi$ is faithful when $n$ is odd, and ${\rm Ker}~\phi$ has order
2 when $n$ is even. From Corollary \ref{c4.2}, $Z(W)\subseteq {\rm
Ker}~\phi,$ and from the fact that a normal subgroup of order 2 is
contained in the center, we have ${\rm Ker}~\phi\subseteq Z(W).$
\end{proof}

\section{Coxeter groups of type $E_n$}\label{s7}
Fix an integer $n\geq 6.$ Let $W$ denote the Coxeter group of type
$E_n$ with the Coxeter graph $S$ in Fig. 2. In this section we shall
determine the orbits of $F_2^n$ under the action of the flipping
group $\mathbf{W}$ of $S.$ Restricting the attention to the case
$n=6, 7$ or $8$ in which $W$ is finite, we show that the kernel of
the Vogan representation $\phi$ of $W$ is the center $Z(W)$ of $W$.
\bigskip

\setlength{\unitlength}{1mm}
\begingroup\makeatletter\ifx\SetFigFont\undefined
\def\x#1#2#3#4#5#6#7\relax{\def\x{#1#2#3#4#5#6}}%
\expandafter\x\fmtname xxxxxx\relax \def\y{splain}
\gdef\SetFigFont#1#2#3{%
  \ifnum #1<17\tiny\else \ifnum #1<20\small\else
  \ifnum #1<24\normalsize\else \ifnum #1<29\large\else
  \ifnum #1<34\Large\else \ifnum #1<41\LARGE\else
     \huge\fi\fi\fi\fi\fi\fi
  \csname #3\endcsname}%
\else \gdef\SetFigFont#1#2#3{\begingroup
  \count@#1\relax \ifnum 25<\count@\count@25\fi
  \def\x{\endgroup\@setsize\SetFigFont{#2pt}}%
  \expandafter\x
    \csname \romannumeral\the\count@ pt\expandafter\endcsname
    \csname @\romannumeral\the\count@ pt\endcsname
  \csname #3\endcsname}%
\fi\endgroup

\begin{picture}(50, 50)
\put(0,25){\footnotesize{$E_n$($n\geq 6$)}}
\put(30,25){\circle{1.5}} \put(38,25){\circle{1.5}}
\put(46,25){\circle{1.5}} \put(54,25){\circle{1.5}}
\put(62,25){\circle{1.5}} \put(46,33){\circle{1.5}}
\put(66,25){\circle*{.8}} \put(70,25){\circle*{.8}}
\put(74,25){\circle*{.8}} \put(78,25){\circle{1.5}}
\put(86,25){\circle{1.5}}\put(94,25){\circle{1.5}}
\put(30.75,25){\line( 1, 0){6.5}}
\put(38.75,25){\line( 1, 0){6.5}}
\put(46.75,25){\line( 1, 0){6.5}}
\put(54.75,25){\line( 1, 0){6.5}}
\put(46,25.75){\line( 0, 1){6.5}}
\put(78.75,25){\line( 1, 0){6.5}}
\put(86.75,25){\line( 1, 0){6.5}}
\put(30,22){\scriptsize{$s_{n-1}$}}
\put(38,22){\scriptsize{$s_{n-2}$}}
\put(46,22){\scriptsize{$s_{n-3}$}}
\put(54,22){\scriptsize{$s_{n-4}$}}
\put(62,22){\scriptsize{$s_{n-5}$}}
\put(46,35){\scriptsize{$s_n$}}
\put(78,22){\scriptsize{$s_3$}}
\put(86,22){\scriptsize{$s_2$}}
\put(94,22){\scriptsize{$s_1$}}

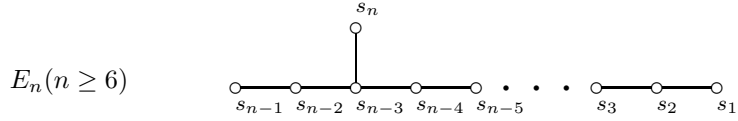
\put(0,5){{\bf Figure 2:} The Coxeter graph of type $E_n.$}
\end{picture}
\bigskip

Set $\overline{1}=\widetilde{s}_1$,
$\overline{i+1}=\mathbf{s_i}\mathbf{s_{i-1}}\cdots
\mathbf{s_1}\overline{1}$ for $1\leq i\leq n-1$ and
$\overline{n+1}=\widetilde{s}_n.$ Note that
\begin{eqnarray}\label{e6.1}
\overline{i}&=&\widetilde{s}_i+\widetilde{s}_{i-1}\quad{\rm for}~
2\leq i\leq n-3,\nonumber\\
\overline{n-2}&=&\widetilde{s}_{n-3}+\widetilde{s}_{n-2}+\widetilde{s}_n,\\
\overline{n-1}&=&\widetilde{s}_{n-2}+\widetilde{s}_{n-1}+\widetilde{s}_n,\nonumber\\
\overline{n}&=&\widetilde{s}_{n-1}+\widetilde{s}_n.\nonumber
\end{eqnarray}
Set $\Delta=\Delta(E_n):=\{\overline{1},
\overline{2},\ldots,\overline{n}\}$ to be the simple basis of
$F_2^n$ in this case. Observe that
\begin{equation}\label{e6.2}
\overline{n+1}=\overline{1}+\overline{2}+\cdots+\overline{n}.
\end{equation}
Set $\Delta(a)$ and $wt(a)$ as before for $a\in F_2^n.$ For example,
$\Delta(\overline{n+1})=\Delta$ and $wt(\overline{n+1})=n.$

\begin{lem}\label{p6.1} The following (i),(ii) hold.
\begin{enumerate}
\item[(i)] For each $1\leq i\leq n-1$,
$\mathbf{s_i}\overline{i}=\overline{i+1}$,
$\mathbf{s_i}\overline{i+1}=\overline{i}$, and
$$\mathbf{s_i}\overline{j}=\overline{j} \qquad {\rm for~~}
\overline{j}\in\{\overline{1},\overline{2},\ldots,\overline{n+1}\}-\{\overline{i},\overline{i+1}\}.$$
\item[(ii)]
$\mathbf{s_n}\overline{n+1}=\overline{n-2}+\overline{n-1}+\overline{n}$,
$\mathbf{s_n}\overline{n}=\overline{n-2}+\overline{n-1}+\overline{n+1}$,
$\mathbf{s_n}\overline{n-1}=\overline{n-2}+\overline{n}+\overline{n+1}$,
$\mathbf{s_n}\overline{n-2}=\overline{n-1}+\overline{n}+\overline{n+1}$
and
$$\mathbf{s_n}\overline{j}=\overline{j} \qquad {\rm for~~} 1\leq
j\leq n-3.$$
\end{enumerate}
\end{lem}
\begin{proof}
This is immediate by applying Lemma \ref{l2.2} and (\ref{e6.1}).
\end{proof}

Let $S_n$ denote the group of permutations on
$\Delta=\{\overline{1},\overline{2},\ldots,\overline{n}\}$. Set
$T:=\{s_1,s_2,\ldots,s_{n-1}\}$. Recall that $\mathbf{W}_T$ is the
subgroup of $\mathbf{W}$ generated by $\{\mathbf{s}~|~s\in T\}$.
By Lemma \ref{p6.1}, we find that the set $\Delta$ is closed under
the left multiplication of elements in $\mathbf{W}_T$.

\begin{defn}\label{d7.2}
Let $\epsilon:\mathbf{W}_T\rightarrow S_n$ denote the homomorphism
satisfying
$$
\epsilon(G)(\overline{j})=G\overline{j}
$$
for $1\leq j\leq n$ and $G\in \mathbf{W}_T.$
\end{defn}

In fact, $\epsilon$ is an isomorphism since $\Delta$ is a spanning
set and the $n-1$ transpositions
$\epsilon(\mathbf{s_1}),\epsilon(\mathbf{s_2}),\ldots,\epsilon(\mathbf{s_{n-1}})$
generate $S_n$.

\begin{prop}\label{p6.2} The following are orbits of $F_2^n$ under the action of
$\mathbf{W}.$
\begin{eqnarray}\label{ee7.3}
O_0&=&\{0\},\nonumber\\
O_1&=&\{a\in F_2^n~|~a\not=0, wt(a)\equiv 1~{\rm or~}
n-2~\pmod{4}\},\\
O_2&=&\{a\in F_2^n~|~a\not=0, wt(a)\equiv 2~{\rm or~}
n-3~\pmod{4}\},\nonumber\\
O_3&=&\{a\in F_2^n~|~a\not=0, wt(a)\equiv 3~{\rm or~}
n~\pmod{4}\},\nonumber\\
O_4&=&\{a\in F_2^n~|~a\not=0, wt(a)\equiv 0~{\rm or~}
n-1~\pmod{4}\}.\nonumber
\end{eqnarray}
In particular $O_1=O_3$ when $n\equiv 1{\pmod 4}$, $O_1=O_4$ and
$O_2=O_3$ when $n\equiv 2{\pmod 4}$, $O_2=O_4$ when $n\equiv
3{\pmod 4}$, and $O_1=O_2$ and $O_3=O_4$ when $n\equiv 0{\pmod
4}$.
\end{prop}
\begin{proof}
It is clear that $O_0$ is an orbit.  There are four cases to put
nonzero vectors $a,$ $b$ in an orbit. (a)$wt(a)=wt(b):$ This is
because
 $\epsilon(\mathbf{W}_T)=S_n$ acts transitively on the fixed size subsets
of $\Delta$; (b) $wt(b)=n+3-wt(a),$ or $n-1-wt(a):$ This is from (a)
and the observation that
\begin{equation}\label{e6.3}
wt(\mathbf{s_n}a)=\left\{
\begin{array}{llr}
n+3-wt(a), & \hbox{if $|\Delta(a)\cap \{\overline{n},\overline{n-1},\overline{n-2}\}|=3$;} \\
n-1-wt(a), & \hbox{if $|\Delta(a)\cap \{\overline{n},\overline{n-1},\overline{n-2}\}|=1$;} \\
w(a), & \hbox{else} \\
\end{array}
\right.
\end{equation}
 by Lemma \ref{p6.1}(ii) and
(\ref{e6.2}); (c)  $wt(a)=wt(b)-4:$ This is by
 applying the first case of
(\ref{e6.3}) and then applying the second case of (\ref{e6.3});
and (d) $wt(a)=wt(b)+4:$ This is by applying the second case of
(\ref{e6.3}) and then the first case of (\ref{e6.3}). The
proposition follows from the above cases (a)-(d).
\end{proof}

\begin{rem}
With reference to Proposition~\ref{p6.2}, for each orbit $O$ of $F_2^n$ with $O\not=O_0$ there is $1\leq i\leq n$ such that $\widetilde{s}_{i}\in O.$ For example $\widetilde{s}_i\in O_i$ for $i=1,2,3$ and $\widetilde{s}_{n-1}\in O_4.$
\end{rem}

Similar to case of $A_n$,  we determine the reducibility of $\phi$
from Proposition \ref{p6.2} immediately.

\begin{prop}
The Vogan representation $\phi$ is irreducible if and only if $n$
is even.\hfill $\Box$
\end{prop}

Recall that for $a\in F_2^n,$ the {\it isotropy group} of $a$ in
$\mathbf{W}$ is $\{G\in \mathbf{W}~|~Ga=a\}$, and the cardinality of
the orbit of $a$ is equal to the index of the isotropy group of $a.$

\begin{cor}\label{c6.1}
For $J:=\{s_2,s_3,\ldots,s_n\}.$ the number $|\mathbf{W}_J||O_1|$
divides $|\mathbf{W}|$, where
\begin{align}\label{e6.4}
|O_1|=\left\{
\begin{array}{ll}
2^{n-1}-(-1)^{\frac{n}{4}}2^{\frac{n-2}{2}}, & \hbox{if $n\equiv 0\pmod{4},$}\\
2^{n-1}, & \hbox{if $n\equiv 1\pmod{4},$ } \\
2^{n-1}+(-1)^{\frac{n-2}{4}} 2^{\frac{n-2}{2}}-1, & \hbox{if $n\equiv 2\pmod{4},$} \\
2^{n-2}+(-1)^{\frac{n-3}{4}} 2^{\frac{n-3}{2}}, & \hbox{if $n\equiv 3\pmod{4}.$} \\
\end{array}
\right.
\end{align}
\end{cor}
\begin{proof}
Since $ \mathbf{W}_J$ is a subgroup of the isotropy group of
$\overline{1}$, the number $|\mathbf{W}_J||O_1|$ divides
$|\mathbf{W}|$. Note that by (\ref{ee7.3})
$$
|O_1|=\left\{
\begin{array}{ll}
\sum\limits_{k\equiv 1,2(\bmod{4})\atop 1\leq k \leq n}{n \choose k}, & \hbox{if $n\equiv 0\pmod{4},$}\\
\sum\limits_{k\equiv 1(\bmod{2})\atop 1\leq k \leq n}{n \choose k}, & \hbox{if $n\equiv 1\pmod{4},$}\\
\sum\limits_{k\equiv 0,1(\bmod{4})\atop 1\leq k \leq n}{n \choose k}, & \hbox{if $n\equiv 2\pmod{4},$}\\
\sum\limits_{k\equiv 1(\bmod{4})\atop 1\leq k \leq n}{n \choose k}, & \hbox{if $n\equiv 3\pmod{4},$}\\
\end{array}
\right.
$$
where $n\choose k$ is the binomial coefficient. From this, we
routinely prove (\ref{e6.4}) by induction on $n$.
\end{proof}

We need to quote a lemma.

\begin{lem}\label{l6.3}(\cite[Lemma 10.2.11]{bcn:89})
If $W$ is of type $E_7$ or $E_8$ then $Z(W)=\{1,w_0\}$, where $w_0$
is the longest element of $W$. \hfill $\Box$
\end{lem}

Recall that $T=\{s_1,s_2,\ldots,s_{n-1}\}$ and
$J=\{s_2,s_3,\ldots,s_n\}$ when we indicate that the Coxeter group
$W$ is of type $E_n.$

\begin{prop}\label{c6.2}
The Vogan representation $\phi$ of $W$ is faithful if $W$ is of
type $E_6$, and  $|{\rm Ker}~\phi|=2$ if $W$ is of type $E_7$.
Moreover, ${\rm Ker}~\phi=Z(W)$ if $W$ is of type $E_6$ or $E_7$.
\end{prop}
\begin{proof}
Suppose $W$ is of type $E_6$. With referring to Corollary
\ref{c6.1}, we have $|O_1|=27$. By Lemma \ref{l3.5}(iii) and
Proposition \ref{p5.12} (the case $D_5$), we know
$|\mathbf{W}_J|=2^45!$, where $J$ is of type $D_5$. Since
$|\mathbf{W}_J||O_1|$ divides $|\mathbf{W}|$, we have
$|\mathbf{W}|\geq 2^45!\cdot27=2^73^45$. Since $|W|=2^73^45$
\cite[p.44]{h:90}, $W$ is isomorphic to $\mathbf{W}$ and ${\rm
Ker}~\phi$ is trivial. By this and Corollary \ref{c4.2},  $Z(W)$ is
trivial.

Suppose $W$ is of type $E_7$. From Corollary \ref{c4.2} and Lemma
\ref{l6.3}, $|{\rm Ker}~\phi|\geq 2$. Since $|W|=2^{10}3^45\cdot7$
\cite[p.44]{h:90}, we see that $|\mathbf{W}|\leq2^93^45\cdot7$. On
the other hand, according to a similar counting argument as above,
we have $|O_1|=28$, $|\mathbf{W}_J|=2^73^45,$ where $J$ is of type
$E_6,$ and hence $|\mathbf{W}|\geq2^93^45\cdot7$. Thus,
$|\mathbf{W}|=2^93^45\cdot7$ and $|Z(W)|=|{\rm Ker}~\phi|=2$.
\end{proof}

We now go to the last case $W$ of type $E_8.$ 
Note that $J=\{s_2, s_3, \ldots, s_8\}$ is of type $E_7$ and
$T\cap J=\{s_2, s_3, \ldots, s_7\}$ is of type $A_6.$ We need more
information of the nontrivial element $w_0$ in the center $Z(W_J)$
of $W_J.$ It is quite complicate to describe $w_0$ directly as a
product of elements in $J.$ We borrow two notations to describe
$w_0.$ Let $\phi$ denote the Vogan representation of $W$. Note
that $\phi\upharpoonright W_{T\cap J}$ is an isomorphism of
$W_{T\cap J}$ onto $\mathbf{W}_{T\cap J}$ by Lemma~\ref{l3.5}(ii)
and Proposition~\ref{p4.5}. Also
$\epsilon\upharpoonright\mathbf{W}_{T\cap J}:\mathbf{W}_{T\cap
J}\rightarrow S_7$ is an isomorphism , where $\epsilon$ is as in
Definition~\ref{d7.2} and $S_7$ is the group of permutations on
$\{\overline{2},\overline{3},\ldots,\overline{8}\}.$ The
expression of $w_0$ is as follows.
\begin{equation}\label{e7.6}
\left.
\begin{array}{ll}
w_0=&\phi^{-1}(\epsilon^{-1}( (\overline{2},
\overline{8},\overline{3}, \overline{7}, \overline{4}, \overline{6},
\overline{5})))s_8\phi^{-1}(\epsilon^{-1}( (\overline{5},
\overline{8})(\overline{4}, \overline{7})(\overline{3},
\overline{6})))s_8\\
&\phi^{-1}(\epsilon^{-1}( (\overline{4}, \overline{8})(\overline{3},
\overline{7})(\overline{2},
\overline{6})))s_8\phi^{-1}(\epsilon^{-1}( (\overline{5},
\overline{8})(\overline{4},
\overline{7})))s_8\\&\phi^{-1}(\epsilon^{-1}( (\overline{3},
\overline{7})(\overline{2}, \overline{6})))s_8.
\end{array}
\right.
\end{equation}
It is routine to check that the above $w_0$ maps to $-I$ by the
faithful representation defined in \cite[Proposition 8]{ab:05}
with $c=0$ or in \cite[p. 291]{he:94} to conclude $w_0$ is in the
center of $W_J$ and indeed is the longest element of $W_J$ by
\cite[Proposition 21]{ab:05}. Thus, we have the following lemma.

\begin{lem}\label{l7.8}
Let $W$ be of type $E_8$ with the Vogan representation $\phi$ and
$w_0\in Z(W_J)$ be not identity. Then $\phi(w_0)$ is
\begin{eqnarray*}
&& \epsilon^{-1}((\overline{2}, \overline{8},\overline{3},
\overline{7}, \overline{4}, \overline{6},
\overline{5}))\mathbf{s_8}\epsilon^{-1}((\overline{5},
\overline{8})(\overline{4}, \overline{7})(\overline{3},
\overline{6}))\mathbf{s_8} \epsilon^{-1}((\overline{4},
\overline{8})(\overline{3}, \overline{7})(\overline{2},
\overline{6}))\mathbf{s_8}\\&\times& \epsilon^{-1}((\overline{5},
\overline{8})(\overline{4},
\overline{7}))\mathbf{s_8}\epsilon^{-1}((\overline{3},
\overline{7})(\overline{2}, \overline{6}))\mathbf{s_8}.
\end{eqnarray*}
\hfill $\Box$
\end{lem}

Note that $W_J$ is not isomorphic to its flipping group
$\mathbf{W_J}$ by Proposition~\ref{c6.2}. The following lemma claims
that $W_J$ is isomorphic to the subgroup $\mathbf{W}_J$ of
$\mathbf{W}$.

\begin{lem}\label{l6.4}
Let $W$ be of type $E_8$ with the Vogan representation $\phi.$
Then the restriction $\phi\upharpoonright W_J$ of $\phi$ to $J$ is
injective.
\end{lem}
\begin{proof}Let
$\phi':W_J\rightarrow \mathbf{W_J}$ denote the Vogan
representation of $W_J$. From Lemma \ref{l3.5}(iii) and
Proposition \ref{c6.2}, we see that ${\rm Ker}~\phi\upharpoonright
W_J\subseteq {\rm Ker}~\phi'=\{1,w_0\}$, where $w_0$ is given in
(\ref{e7.6}). To prove that ${\rm Ker}~\phi\upharpoonright W_J$ is
trivial, it suffices to show that $\phi(w_0)\not=I.$ This follows
from the computation
$$
\phi(w_0)\overline{8}=\overline{1}+\overline{8}
$$
by applying the expression $\phi(w_0)$ in Lemma~\ref{l7.8} to
$\overline{8}$ and using  Lemma~\ref{p6.1} and (\ref{e6.2}) for
$n=8$ to simplify.
\end{proof}

There is a similar result about $W$ of type $E_8.$

\begin{prop}
If $W$ is of type $E_8$ then ${\rm Ker}~\phi$ has order 2. Moreover,
${\rm Ker}~\phi=Z(W)$.
\end{prop}
\begin{proof}
We have $|O_1|=2^3\cdot 3\cdot 5$ from (\ref{e6.4}),
$|\mathbf{W}_J|=|W_J|=2^{10}3^45\cdot7$ from  Lemma \ref{l6.4} and
$|W|=2^{14}3^55^27$\cite[p.44]{h:90}. Therefore, as the proof of
Proposition \ref{c6.2},  ${\rm Ker}~\phi$ has order 2 and ${\rm
Ker}~\phi=Z(W)$.
\end{proof}

\section{Concluding remarks}

We list the main results of this article as follows.
\bigskip

\begin{center}
\begin{tabular}{ccc}
\hline Dynkin diagram   &reducibility of $\phi$ &$|{\rm Ker}~\phi|$  \\
\hline \\[-8pt]
$A_n$ &{\small $\phi$ is irr. iff $n=1$ or $n$ is even.} &{\small
$\left\{
\begin{array}{ll}
2,& \hbox{if $n=1,$}\\
1,& \hbox{else.} \\
\end{array}
\right. $ }
\\
\hline \\[-8pt]
$\left.
\begin{array}{c}
D_n\\
 (n\geq 4)\\
\end{array}
\right.$
 &{\small $\phi$ is not irr.} &{\small$\left\{
\begin{array}{ll}
2,& \hbox{if $n$ is even,}\\
1,& \hbox{else.} \\
\end{array}
\right.$}
\\
\hline \\[-8pt] $E_6$  &{\small $\phi$ is irr.} &$1$
\\

\hline \\[-8pt] $E_7$  &{\small $\phi$ is not irr.} &$2$
\\

\hline \\[-8pt] $E_8$  &{\small $\phi$ is irr.} &$2$
\\
\hline\\[3pt]
\end{tabular}
{\bf Table 1:} The reducibility and the kernel of a Vogan
representation $\phi.$
\end{center}
 \bigskip

\begin{center}
\begin{tabular}{cc}
\hline Coxeter graph  &orbits\\
\hline \\[-8pt] $A_n$  &{\small$
\left.
\begin{array}{l}
O_i=\{a\in F_2^n~|~wt(a)=i~{\rm or}~n+1-i\}
\hbox{$(0 \leq i\leq \lfloor\frac{n+1}{2}\rfloor$).}\\
\end{array}
\right. $}
\\
\hline \\[-8pt]
$\left.
\begin{array}{c}
D_n\\
 (n\geq 4)\\
\end{array}
\right.$
&{\small $ \left.
\begin{array}{l}
O_i=\{a\in Z~|~wt(a)=i~{\rm or}~n-i\}~~~(0\leq i\leq \lfloor\frac{n}{2}\rfloor),\\
\Omega_o=\{a\in F_2^n-Z~|~wt(a)\equiv 1~{\rm
or}~n-1~\pmod{2}\},\\
\Omega_e=\{a\in F_2^n-Z~|~wt(a)\equiv 0~{\rm or}~n~\pmod{2}\},\\
\hbox {$\Omega_o=\Omega_e=F_2^n-Z$ when $n$ is odd.}\\
\end{array}
\right. $}
\\
\hline \\[-8pt]
$\left.
\begin{array}{c}
E_n\\
 (n\geq 6)\\
\end{array}
\right.$
 &{\small $ \left.
\begin{array}{l}
O_0=\{0\},\\
 O_1=\{a\in F_2^n~|~a\not=0, wt(a)\equiv 1~{\rm or~}
n-2~\pmod{4}\},\\
O_2=\{a\in F_2^n~|~a\not=0, wt(a)\equiv 2~{\rm or~}
n-3~\pmod{4}\},\\
O_3=\{a\in F_2^n~|~a\not=0, wt(a)\equiv 3~{\rm or~}
n~\pmod{4}\},\\
O_4=\{a\in F_2^n~|~a\not=0, wt(a)\equiv 0~{\rm or~}
n-1~\pmod{4}\}.\\
\hbox{$O_1=O_3$ when $n\equiv 1{\pmod
4}$,}\\
\hbox{$O_1=O_4$ and $O_2=O_3$ when $n\equiv 2{\pmod 4}$,}\\
\hbox{$O_2=O_4$ when $n\equiv 3{\pmod 4}$,}\\
\hbox{$O_1=O_2$ and $O_3=O_4$ when $n\equiv 0{\pmod 4}$}.\\
\end{array}
\right.$}
\\

\hline \\[3pt]
\end{tabular}
{\bf Table 2:} The orbits of $F_2^n$ under the action of the
flipping group of a Coxeter graph $S.$
\end{center}

\bigskip

\noindent Hau-wen Huang \hfil\break Department of Applied
Mathematics \hfil\break National Chiao Tung University \hfil\break
1001 Ta Hsueh Road \hfil\break Hsinchu, Taiwan 30050,
R.O.C.\hfil\break Email: {\tt poker80@msn.com} \hfil\break Fax:
+886-3-5724679 \hfil\break
\medskip

\end{document}